\newif\iftechreport
\techreporttrue

\documentclass[12pt]{l4dc2021} 
\usepackage{amsmath}
\usepackage{amssymb}
\usepackage{csquotes}
\usepackage{tikz}
\usetikzlibrary{automata, positioning, arrows, patterns}
\usepackage{pgfplots}
\pgfplotsset{compat=1.11} 

\usepackage[capitalize]{cleveref}

\newtheorem{assumption}{Assumption}
\newtheorem{problem}{Problem}

\newcommand{\BnB}{Branch and bound}
\newcommand{\bnb}{branch and bound}

\newcommand{\Hsys}{\mathcal{S}}
\newcommand{\hsyst}{discrete-time control hybrid system}
\newcommand{\Hsyst}{Discrete-time control hybrid system}
\newcommand{\states}{states}

\newcommand{\lenub}{L}

\title{Abstraction-based branch and bound approach to Q-learning for hybrid optimal control}
\usepackage{times}

\author{%
 \Name{Beno\^it Legat} \Email{benoit.legat@uclouvain.be}\\
 \Name{Rapha\"el M. Jungers} \Email{raphael.jungers@uclouvain.be}\\
 \addr ICTEAM, UCLouvain, 4 Av. G. Lema\^itre, 1348 Louvain-la-Neuve, Belgium
 \AND
 \Name{Jean Bouchat} \Email{jean.bouchat@uclouvain.be}\\
 \addr ELI, UCLouvain, 2 Croix du Sud, 1348 Louvain-la-Neuve, Belgium
}

\newcommand{\defref}[1]{Definition~\ref{#1}}
\newcommand{\propref}[1]{Proposition~\ref{prop:#1}}
\newcommand{\lemref}[1]{Lemma~\ref{lem:#1}}

\newcommand{\modename}{mode}
\newcommand{\dist}{d}
\newcommand{\distx}{d'}
\newcommand{\modes}{\mathcal{Q}}
\newcommand{\mode}{q}

\newcommand{\modeT}{q_t}
\newcommand{\statespace}{\mathcal{X}}
\newcommand{\dynto}{\rightsquigarrow}
\newcommand{\dyntou}[1]{\stackrel{#1}{\rightsquigarrow}}
\newcommand{\horizon}{H}
\newcommand{\deadline}{T}
\newcommand{\cost}{k}

\newcommand{\QcV}[1]{\mathcal{T}^Q_{c}{#1}}
\newcommand{\Qfun}{Q}
\newcommand{\Qfunv}{\hat{Q}}
\newcommand{\Qv}[3]{\Qfunv((#1, #2), #3)}
\newcommand{\Qvu}[5][\Qfun]{{#1}((#2, #3), (#4, #5))}
\newcommand{\policy}{\mu}

\DeclareMathOperator{\Dom}{Dom}
\newcommand{\valuefun}{J}

\DeclareMathOperator*{\argmax}{arg\,max}
\DeclareMathOperator*{\argmin}{arg\,min}
\DeclareMathOperator*{\Post}{Post}

\usepackage{ifthen}
\newcommand{\emptytuple}{\emptyset}
\newcommand{\tuplecat}[2]{(#1; #2)}
\newcommand{\tuple}[3]{\ifthenelse{\equal{#2}{1}}{\mathbf{#1}_{#3}}{({#1}_i)_{i=#2}^{#3}}}
\newcommand{\tuplep}[3]{\ifthenelse{\equal{#2}{1}}{\mathbf{#1}'_{#3}}{({#1}'_i)_{i=#2}^{#3}}}
\newcommand{\tuples}[3]{\ifthenelse{\equal{#2}{1}}{\mathbf{#1}^\star_{#3}}{({#1}^\star_i)_{i=#2}^{#3}}}

\usepackage{comment}
\usepackage{environ}

\usepackage{xcolor}
\definecolor{troyes}{rgb}{0.996078431372549, 0.9921568627450981, 0.9411764705882353}
\definecolor{frambo}{rgb}{0.7803921568627451, 0.17254901960784313, 0.2823529411764706}
\definecolor{lichen}{rgb}{0.5215686274509804, 0.7568627450980392, 0.49411764705882355}
\definecolor{canard}{rgb}{0.01568627450980392, 0.5450980392156862, 0.6039215686274509}
\definecolor{aurore}{rgb}{1.0, 0.796078431372549, 0.3764705882352941}
\definecolor{re}{rgb}{0.9607843137254902, 0.5725490196078431, 0.592156862745098}
\definecolor{gre}{rgb}{0.796078431372549, 0.8745098039215686, 0.5019607843137255}
\definecolor{blu}{rgb}{0.13725490196078433, 0.5803921568627451, 0.807843137254902}
\definecolor{ora}{rgb}{0.9803921568627451, 0.796078431372549, 0.5843137254901961}
\definecolor{yel}{rgb}{0.9490196078431372, 0.9411764705882353, 0.5450980392156862}

\iftechreport%
\newenvironment{myproof}{\begin{proof}}{\end{proof}}
\else%
\NewEnviron{myproof}{}{}
\fi%

\newenvironment{figol}
{
\begin{tikzpicture}[thick]
    \begin{axis}[
        height=0.48\textwidth,
        width=0.48\textwidth,
        xmin=-10,
        xmax=2.5,
        xlabel={$x_1$},
        ymin=-10,
        ymax=2.5,
        ylabel={$x_2$},
        axis lines=left,
        every axis x label/.style={at={(ticklabel* cs:1)},anchor=west},
        every axis y label/.style={at={(ticklabel* cs:1)},anchor=south},
    ]    
        \draw[dashed] (1.85,-10)--(1.85,2);
        \draw[dashed] (-10,2)--(1.85,2);
        \draw (-0.5,-0.5) rectangle (0.5,0.5) node[pos=0.5] {$\mathcal{X}_t$};
        \draw [pattern=north west lines, pattern color=gray] (-10,-10) rectangle (-5,-5) node[pos=0.5] {$\mathcal{O}_1$};
        \draw [pattern=north west lines, pattern color=gray] (-5,-4) rectangle (1.85,-3) node[pos=0.5] {$\mathcal{O}_2$};
        \draw (-6,1) rectangle (-5,2) node[pos=0.5] {$\mathcal{A}$};
        \draw (-5,-3) rectangle (-4,-2) node[pos=0.5] {$\mathcal{B}$};
}
{
    \end{axis}
\end{tikzpicture}
}

\usepackage{siunitx}

\begin{document}
\maketitle

\begin{abstract}%

In this paper, we design a theoretical framework allowing to apply model predictive control on hybrid systems.
For this, we develop a theory of approximate dynamic programming
by leveraging the concept of alternating simulation.
We show how to combine these notions in a branch and bound algorithm
that can further refine the Q-functions using Lagrangian duality.
We illustrate the approach on a numerical example.

\end{abstract}


\section{Introduction}

The capability of hybrid systems to model both continuous dynamics and discrete events in the same mathematical model
renders them essential in fields such as robotics, automotive control or air traffic management.
However, with their ability to model such complex systems come substantial challenges for controlling them.
In this work\footnote{
  A version of this paper containing the proofs is available in~\cite{legat2020abstractionCO}.
  The results of the numerical experiments presented in \cref{sec:example} can be reproduced
  using the Code Ocean capsule in \cite{legat2020abstractionarXiv}.
  It relies on the OSQP solver~\cite{osqp} for solving quadratic programs through the MathOptInterface~\cite{legat2020mathoptinterface}.
}, we study finite time horizon optimal control problems on hybrid systems.

For a linear hybrid system, a quadratic cost function and a fixed choice of discrete control inputs,
the optimal value of the continuous control inputs can be found solving a Quadratic Program (QP) \cite{bemporad2002explicit}.
However, the number of discrete control inputs typically grows exponentially with the time horizon or ``size'' of the system.
In \cite{bemporad1999control}, the authors introduce a Mixed Integer Quadratic Program (MIQP) that simultaneously
finds the optimal value of both the discrete and continuous control inputs.
While the number of integer variables of the MIQP grows linearly with the time horizon or ``size'' of the system,
MIQPs are NP-hard to solve in general hence this approach 
is not suitable for the online control of large-sized problems with real-time constraints.

Several approaches were proposed to enable a small horizon Model Predictive Controller (MPC) to
satisfy such real-time constraints online along with the control objective.
In \cite{gol2013language, gol2015temporal}, the authors develop an algorithm to obtain a Lyapunov function
that guarantees the MPC controller to reach a target discrete state.
Computing this Lyapunov function can however be prohibitive and their method is not ensured to find
an optimal solution. 
In \cite{bouchatreinforcement, menta2020learning}, the authors show how the weak duality of the MIPQ allows
the refinement an under-approximation given by a value or Q-function.
This value or Q-function can be used as terminal cost of the MPC to improve the cost of the solution found.

Computing a Lyapunov function or an accurate approximation of the value or Q-function over the whole state-space is intractable for most classes of hybrid systems \cite{blondel2000survey},
hence it seems appropriate to only generate an accurate approximation along the optimal trajectory.
As the optimal trajectory is unknown, \cite{bouchatreinforcement} alternates between
1) a search for a sub-optimal trajectory according to the current under-approximation of the value function
using Model Predictive Control (MPC), which they called \emph{forward pass}, and
2) a refinement of the under-approximation of the value function along the trajectory, called \emph{backward pass}.

In \cref{sec:abstraction}, we formalize an approach based on \emph{simulation} relations to obtain \emph{Lyapunov} functions and \emph{Bellman-like} Q-functions.
This generalizes the algorithm of \cite{gol2013language} for Lyapunov functions. 
This abstraction approach provides both a Bellman-like value function on the whole state-space as well as a Lyapunov function
on some 
set $\mathcal{X}_f$ containing the target set.

In \cref{sec:mip}, we show how to combine a Lyapunov and a Bellman-like Q-function in a \bnb{} algorithm
solving an optimal control problem.
Since learning an accurate Q-function in the whole state-space is intractable,
the algorithm only refines it along trajectories computed with an MPC-approach throughout the algorithm.

In \cref{sec:example}, we demonstrate the algorithm on an example from \cite{gol2013language, bouchatreinforcement} illustrated in \cref{fig:plots}.
\begin{figure}[htbp]
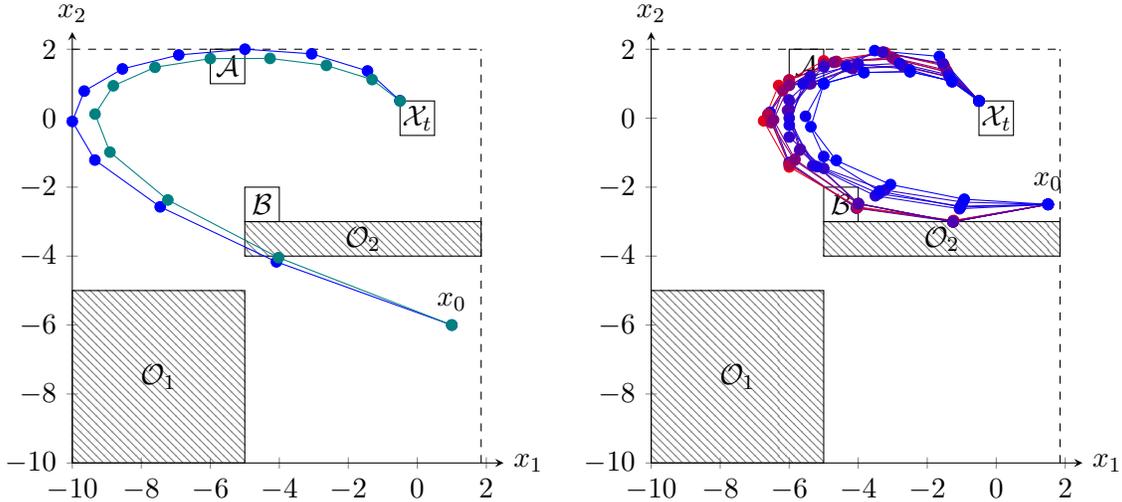

\floatconts
  {fig:plots}
  {\caption{Feasible trajectories providing increasingly better upper bounds found by \cref{algo:BnB}
  on the example detailed in \cref{sec:example}.
  The left (resp. right) figure provides the trajectories found for instance $\mathcal{I}_1$ (resp. $\mathcal{I}_2$) with Bellman-like Q-function $Q_2$;
  see \cref{sec:example} for the definition of $\mathcal{I}_1, \mathcal{I}_2$ and $Q_2$.}}
  {%
\begin{figol}
        \addplot[mark=*,color=blue] coordinates
		{ (1.0,-6.0) (-4.08326083652693, -4.16652157792225) (-7.452103692575193, -2.571163630780764) (-9.344647629427765, -1.2139233592305483) (-9.999004371331917, -0.09478891152710031) (-9.653269015312578, 0.7862610956032603) (-8.545511305383734, 1.4292559716987285) (-6.913769471419556, 1.8342294340939918) (-4.996048425321604, 2.00121442419686) (-3.0616913911649153, 1.867491079841009) (-1.4429744774058832, 1.3699304105343828) (-0.5037481670764479, 0.5085061781102557) }; 
		\addplot[mark=*,color=teal] coordinates
		{ (1.0,-6.0) (-4.022047002972235, -4.044093992822277) (-7.230161685513392, -2.372135358355666) (-8.90829115789696, -0.9841235768359475) (-9.340382570533198, 0.11994075050249883) (-8.810385426650052, 0.9400535230952529) (-7.602256637656267, 1.4762040350523633) (-5.999968566683014, 1.728372102463325) (-4.26866193347106, 1.7342410887440363) (-2.635795880416679, 1.5314909667308876) (-1.3100047622834396, 1.1200912560675875) (-0.4999518407560096, 0.5000146412394229) }; 
		\node[circle,fill=black,inner sep=0pt,minimum size=3pt,label=above:{$x_0$}] (a) at (1,-6) {};
\end{figol}~%
\begin{figol}
\addplot[mark=*,color=blue] coordinates
{ (1.5, -2.5) (-1.2499358754855439, -2.999873160529467) (-4.053813913222866, -2.6078861631133274) (-6.001707356443967, -1.287900546490158) (-6.561521739419388, 0.16827593846507938) (-5.992898756779733, 0.9689742102170001) (-5.005728822619895, 1.0053491453370698) (-3.5229529469125835, 1.9602022184972403) (-1.6468883584070424, 1.7919258045004578) (-0.5006107347288019, 0.5006276324901333) };
\addplot[mark=*,color=blue!10!red] coordinates
{ (1.5, -2.5) (-1.2337120828649804, -2.9674270318238234) (-4.003850916858565, -2.5728548357681653) (-5.99701190170787, -1.4134611030961888) (-6.7419116691231356, -0.07633883305219247) (-6.305759724619973, 0.9486418662365252) (-5.000654493562242, 1.6615673000716265) (-3.217992987946002, 1.9037559434250266) (-1.507873580944922, 1.5164833502475041) (-0.49976939612242394, 0.49972574910687134) };
\addplot[mark=*,color=blue!20!red] coordinates
{ (1.5, -2.5) (-1.2407445153637375, -2.9814891616572305) (-4.029697141055694, -2.5964161269016777) (-6.000296124707743, -1.3447824899063459) (-6.617422121562451, 0.11053063364717125) (-6.0088294787042305, 1.1066548475244733) (-4.633690199904898, 1.6436238507600083) (-2.9511363002970206, 1.7214839966200797) (-1.42025251551708, 1.3402835106688569) (-0.5000783387285871, 0.5000646485198628) };
\addplot[mark=*,color=blue!30!red] coordinates
{ (1.5, -2.5) (-1.2410842067243961, -2.9821552955104) (-4.029965847220101, -2.595576562422261) (-5.997904800142269, -1.3403106531914148) (-6.609608964068851, 0.11689151074682466) (-5.99660275944713, 1.1091197462945162) (-4.621581599618624, 1.6409166056499251) (-2.9427464546104503, 1.7167512951219492) (-1.4162960420332653, 1.3361500591952145) (-0.4987873574706015, 0.4988703918010905) };
\addplot[mark=*,color=blue!40!red] coordinates
{ (1.5, -2.5) (-1.2527968928142417, -3.005583131348277) (-3.9955259211030705, -2.4798656897663953) (-5.847138263204088, -1.223359011311763) (-6.523547954848565, -0.12946071351305471) (-6.187237490127514, 0.8020809746678754) (-5.000497148006229, 1.5713987274492056) (-3.259219474367183, 1.9111568232902085) (-1.5266622868764899, 1.553957967477596) (-0.4997999200925019, 0.49976739573695905) };
\addplot[mark=*,color=blue!50!red] coordinates
{ (1.5, -2.5) (-1.2524056418886604, -3.004801824087647) (-3.9958777353374164, -2.4821339962880016) (-5.8329119961407905, -1.1919348142972992) (-6.452590827006608, -0.04742356214349455) (-6.000536957829477, 0.9515301584848129) (-4.71578830063719, 1.6179672025865202) (-3.024375740398228, 1.7648580847418793) (-1.44587105859569, 1.3921515673461324) (-0.4998808096909799, 0.49982933954265985) };
\addplot[mark=*,color=blue!60!red] coordinates
{ (1.5, -2.5) (-1.2534013677517586, -3.0067903368272093) (-3.9949187457062854, -2.4762339895299954) (-5.697817750350962, -0.9295640022194334) (-6.047972423301555, 0.2292543635557122) (-5.432975149385742, 1.0007395396234673) (-4.205221213990454, 1.4547680432648082) (-2.7165518467542897, 1.5225699981579175) (-1.3530026821452041, 1.2045272364359056) (-0.5003690128323959, 0.500738624629213) };
\addplot[mark=*,color=blue!70!red] coordinates
{ (1.5, -2.5) (-1.253122027190447, -3.0062310770412823) (-3.993838204228826, -2.4751894127833527) (-5.679534729872405, -0.8962051384446325) (-6.001089431809266, 0.25309317625852773) (-5.374374434660684, 1.0003364760813374) (-4.155509906443693, 1.4373926229106524) (-2.686963119024582, 1.4997010446816428) (-1.343508157695244, 1.187209021081415) (-0.4999526077697524, 0.4999022710831247) };
\addplot[mark=*,color=blue!80!red] coordinates
{ (1.5, -2.5) (-0.9675487743698932, -2.435096792441424) (-3.228356161018757, -2.086516806563938) (-4.998715348328627, -1.4541999988423047) (-5.998797334263548, -0.5459658657960521) (-6.005886776009843, 0.5317805957496473) (-4.995291549434857, 1.4894190587969456) (-3.293450173945107, 1.914262565258025) (-1.5439085969460684, 1.584818512849521) (-0.500912617955579, 0.501170440890695) };
\addplot[mark=*,color=blue!85!red] coordinates
{ (1.5, -2.5) (-1.0633413076497993, -2.6266826094346847) (-3.5024020295170373, -2.251438822172047) (-5.315254286729027, -1.3742656697490534) (-5.999965928871841, 0.004842431526609958) (-5.39172684556087, 1.211635707729704) (-3.999981443505376, 1.5718551105186978) (-2.508700246107624, 1.4107072832493919) (-1.276645560819619, 1.0534021427230837) (-0.4999732150096609, 0.499942654784349) };
\addplot[mark=*,color=blue!90!red] coordinates
{ (1.5, -2.5) (-1.026758981148554, -2.5535181185797082) (-3.397661210541276, -2.188286576936407) (-5.194048613019767, -1.4044867413391722) (-5.997686113004784, -0.20278682090670447) (-5.5995638228118985, 0.999028544600066) (-4.347985692494822, 1.504125211780852) (-2.8008168179739803, 1.59021091780268) (-1.377625710986096, 1.2561703614049091) (-0.4990854015950138, 0.5009106197334325) };
\addplot[mark=*,color=blue!95!red] coordinates
{ (1.5, -2.5) (-1.0392959134742605, -2.5785916446234696) (-3.3863815503117785, -2.1155793479153537) (-4.999639682037996, -1.1109365484264633) (-5.527783540376893, 0.05464820161001361) (-5.0002380021166495, 1.0004419772123918) (-3.837461391470813, 1.3251112429761849) (-2.4999413372362698, 1.3499288280102917) (-1.287522305524183, 1.074909150464486) (-0.5000398439736347, 0.5000556390883238) };
\addplot[mark=*,color=blue] coordinates
{ (1.5, -2.5) (-0.924952392534689, -2.3499052061955625) (-3.0621531133459996, -1.924497036124217) (-4.636322716887948, -1.223842900877931) (-5.37242918178286, -0.2483693841286432) (-4.996116282509809, 1.0009977902023341) (-3.83327454697361, 1.3246808037431368) (-2.496176913494163, 1.3495127637822153) (-1.2843356618552313, 1.074171107195737) (-0.4983962459000687, 0.49771260959395675) };
\node[circle,fill=black,inner sep=0pt,minimum size=3pt,label=above:{$x_0$}] (a) at (1.5,-2.5) {};
\end{figol}

  }
\end{figure}

\section{Discrete optimal control}
\label{sec:abstraction}

In this section, we define \emph{simulation relations} between discrete-time systems and show how to deduce
a \emph{Lyapunov} function for a system from a Lyapunov function for a simulated system as well
as a \emph{Bellman-like} value function for a system from a Bellman-like value function for a simulation.


We use the following notation for discrete-time control systems.
\begin{definition}[Discrete-time control system]\label{def:discrete_time}
    A \emph{discrete-time control system} is defined as a triple
    $S = (\mathcal{X}, \mathcal{U}, \dynto)$
    where $\mathcal{U}$ is the set of input sets
    and $\dynto$ is the subset of transitions $(x, u, x')$ such that the system can reach $x' \in \mathcal{X}$ from $x \in \mathcal{X}$ with input $u \in \mathcal{U}$.
    We denote $(x, u, x') \in \ \dynto$ as $x \dyntou{u} x'$, and the set of $x'$ such that $x \dyntou{u} x'$ as $\Post_u(x)$.
\end{definition}

We say that a \hsyst{} is \emph{deterministic} if for every state $x \in \mathcal{X}$ and control input $u \in \mathcal{U}$, $\Post_u(x)$ is either empty or a singleton.

The simulation used in this section is commonly referred to as an \emph{alternating simulation}.
\begin{definition}[{Alternating simulation relation~\cite[Definition~4.19 and Definition~4.22]{tabuada2009verification}}]
    \label{def:altsim}
    Consider discrete-time control systems $S_1 = (\mathcal{X}_1, \mathcal{U}_1, \dynto_1)$ and $S_2 = (\mathcal{X}_2, \mathcal{U}_2, \dynto_2)$, as defined in \defref{def:discrete_time}.
    Given a relation $R \subseteq \mathcal{X}_1, \mathcal{X}_2$, consider the \emph{extended relation}
    $R^e$ defined by the set of $(x_1, x_2, u_1, u_2)$ such that
    for every $x_2' \in \Post_{2,u_2}(x_2)$,
    there exists $x_1' \in \Post_{1,u_1}(x_1)$ such that $(x_1', x_2') \in R$.
    If for all $(x_1, x_2) \in R$, and for all $u_1 \in \mathcal{U}$,
    there exists $u_2 \in \mathcal{U}$ such that $(x_1, x_2, u_1, u_2) \in R^e$
    then $R$ is an \emph{alternating simulation relation} and $R^e$ is its associated
    \emph{extended alternating simulation relation}.
\end{definition}

\subsection{Bellman-like value and Q-functions}
In this section, we define \emph{Bellman-like value functions} and \emph{Bellman-like Q-functions}, and show how a Bellman-like value function
of a system can be deduced from the Bellman-like value function of an alternating simulation.
The Bellman-like value function will be used to provide lower bounds for the \bnb{} algorithm in \cref{sec:mip}.

We denote the empty tuple as $\emptytuple{}$,
the $l$-tuple $(u_i)_{i=1}^l$ as $\tuple{u}{1}{l}$
and the concatenation of tuples $\tuple{u}{1}{l_1}, \tuplep{u}{1}{l_2}$ as the $(l_1+l_2)$-tuple $\tuplecat{u}{u'}$.
A \emph{cost function} is a given function $c : \mathcal{X} \times \mathcal{U} \to \mathbb{R}$ and
a \emph{value function} is a function $V : \mathcal{X} \to \mathbb{R}$.
Given a \emph{Q-function} $Q : \mathcal{X} \times \mathcal{U}^l \to \mathbb{R}$ for some $l \in \mathbb{N}$ with some cost function $c$, we recursively define the value of $Q(x, \tuple{u}{1}{l'})$ for $l' > l$ with the following
identity for $k = l+1, \ldots, l'$:
\begin{equation}
   \label{eq:Qrec}
   Q(x, \tuple{u}{1}{k}) = c(x, u_1) + \max_{x' \in \Post_{u_1}(x)} Q(x', \tuple{u}{2}{k}).
\end{equation}
Given a cost function $c$ and a value function $V$, $\QcV{V}$ denotes the Q-function with cost function $c$ such that $\QcV{V}(x, \emptytuple{}) = V(x)$. 

The \emph{Bellman operator} $\mathcal{T}_c$ is defined as\footnote{Note that we have $\mathcal{T}_c V(x) = \mathcal{T}_c \QcV{V}(x, \emptytuple{})$.}:
\begin{align}
  \label{eq:Tc}
  \mathcal{T}_c V(x) & = \min_{u \in \mathcal{U}} \QcV{V}(x, u) &
  \mathcal{T}_c Q(x, \tuple{u}{1}{l}) & = \min_{u' \in \mathcal{U}} Q(x, \tuplecat{\tuple{u}{1}{l}}{u'}).
\end{align}

\begin{definition}[Bellman-like value function]
    \label{def:Bellman-like}
    Consider a discrete-time control system $S = (\mathcal{X}, \mathcal{U}, \dynto)$.
    A value function $V$ is a \emph{Bellman-like value function} of $S$ with cost function $c$
    if $V(x) \le \mathcal{T}_c V(x)$ for all $x \in \mathcal{X}$.
\end{definition}


\begin{definition}[Bellman-like Q-function]
    \label{def:Bellman-likeQ}
    Consider a discrete-time control system $S = (\mathcal{X}, \mathcal{U}, \dynto)$.
    A function $Q$ is a \emph{Bellman-like Q-function} of $S$ with cost function $c$ if
    $Q(x, \tuple{u}{1}{k}) \le Q(x, \tuple{u}{1}{l})$ for all $k, l \in \mathbb{N}, x \in \mathcal{X}, \tuple{u}{1}{l} \in \mathcal{U}^{l}$ such that $k \le l$.
\end{definition}

\begin{proposition}
   \label{prop:Bellman-likeQ}
   Consider a discrete-time control system $S = (\mathcal{X}, \mathcal{U}, \dynto)$.
   If $V$ is a Bellman-like value function of $S$ with cost function $c$,
   then $\QcV{V}$ is a Bellman-like Q-function of $S$ with cost function $c$.
   \begin{myproof}
      Consider $l \in \mathbb{N}, x \in \mathcal{X}, u \in \mathcal{U}$ and $\tuplep{u}{1}{l} \in \mathcal{U}^l$.
      For any $x' \in \Post_{u}(x)$,
      \defref{def:Bellman-like} implies that $V(x') \le \mathcal{T}_c^l V(x')$
      and \eqref{eq:Tc} implies that $\mathcal{T}_c^l V(x') \le \QcV{V}(x', \tuplep{u}{1}{l})$.
      Therefore, we have
      \[
          \QcV{V}(x, u)
          = c(x, u) + \max_{x' \in \Post_{u}(x)} V(x')
           \le c(x, u) + \max_{x' \in \Post_{u}(x)} \QcV{V}(x', \tuplep{u}{1}{l})
           = \QcV{V}(x, \tuplecat{u}{\tuplep{u}{1}{l}}).
      \]
   \end{myproof}
\end{proposition}

\begin{theorem}
    \label{theo:bellman}
    Consider discrete-time control systems $S_1 = (\mathcal{X}_1, \mathcal{U}_1, \dynto_1)$, $S_2 = (\mathcal{X}_2, \mathcal{U}_2, \dynto_2)$, as defined in \defref{def:discrete_time},
    and an alternating simulation relation $R$ such that for each $x_1 \in \mathcal{X}_1$, there is exactly
    one $x_2 \in \mathcal{X}_2$ such that $(x_1, x_2) \in R$, which we denote by $R(x_1)$.
    Given a cost function $c_1$ for $S_1$,
    consider an associated cost function satisfying
    \begin{equation}
        \label{eq:cost}
        c_2(x_2, u_2) \le \min_{(x_1, x_2, u_1, u_2) \in R^e} c_1(x_1, u_1).
    \end{equation}
    If $V_2(x)$ is a Bellman-like value function for $S_2$ with cost function $c_2$,
    then $V_1(x_1) = V_2(R(x_1))$ is a Bellman-like value function for $S_1$ with cost function $c_1$.
    \begin{myproof}
        By \defref{def:altsim}, for all $x_1 \in \mathcal{X}_1$ and $u_1 \in \mathcal{U}_1$, there exists $u_2 \in \mathcal{U}_2$ such that
        $(R(x_2), x_2, u_1, u_2) \in R^e$.
        The inequality in \defref{def:Bellman-like} ensures that there exists $x_2' \in \Post_{2, u_2}(R(x_1))$ such that
        \begin{equation}
            \label{eq:Bellman-like_proof}
            V_2(R(x_1)) \le V_2(x_2') + c_2(R(x_1), u_2).
        \end{equation}
        As $(R(x_2), x_2, u_1, u_2) \in R^e$, there exists $x_1' \in \Post_{u_1}(x_1)$ such that $R(x_1') = x_2'$.
        By \eqref{eq:Bellman-like_proof} and \eqref{eq:cost}, we have
        \begin{align*}
            V_1(x_1) = V_2(R(x_1))
            & \le V_2(R(x_1')) + c_2(R(x_1), u_2)\\
            & \le V_1(x_1') + c_1(x_1, u_1).
        \end{align*}
    \end{myproof}
\end{theorem}

\subsection{Lyapunov functions and receding horizon control}

In this section, we define \emph{Lyapunov functions} and show how a Lyapunov function
of a system can be deduced from a Lyapunov function of an alternatingly simulated system.
We then show how a Lyapunov function can ensure that a model predictive controller or receding horizon controller reaches a target.

\begin{definition}[Lyapunov function]
    \label{def:lyap}
    Consider a discrete-time control system $S = (\mathcal{X}, \mathcal{U}, \dynto)$,
    a set $\mathcal{X}_f \subseteq \mathcal{X}$ and
    a cost function $c$.
    A value function $L$ is a \emph{Lyapunov function} with cost function $c$ for $S$ in $\mathcal{X}_f$ if,
    for all $x \in \mathcal{X} \setminus \mathcal{X}_f$, $L(x) = \infty$, and
    for all $x \in \mathcal{X}_f$, $L(x)$ is finite and $L(x) \ge \mathcal{T}_c L(x)$.
\end{definition}

\begin{theorem}
    Consider discrete-time control systems $S_1 = (\mathcal{X}_1, \mathcal{U}_1, \dynto_1)$, $S_2 = (\mathcal{X}_2, \mathcal{U}_2, \dynto_2)$, as defined in \defref{def:discrete_time},
    and an alternating simulation relation $R$ such that for each $x_2 \in \mathcal{X}_2$, there is exactly
    one $x_1 \in \mathcal{X}_1$ such that $(x_1, x_2) \in R$, which we denote by $R(x_2)$.
    Given a cost function $c_2$ for $S_2$,
    consider an associated cost function satisfying
    \begin{equation}
        \label{eq:decrease}
        c_1(x_1, u_1) \ge \max_{(x_1, x_2, u_1, u_2) \in R^e} c_2(x_2, u_2).
    \end{equation}
    If $L_1(x)$ is a Lyapunov function for $S_1$ with cost function $c_1$,
    then $L_2(x_2) = L_1(R(x_2))$ is a Lyapunov function for $S_2$ with cost function $c_2$.
    \begin{myproof}
        Given $x_2 \in \mathcal{X}_2$, \defref{def:lyap} ensures the existence of $u_1 \in \mathcal{U}_{1}$
        such that
        \begin{equation}
            \label{eq:lyap_proof}
            L_1(x_1') + c_1(R(x_2), u_1) \le L_1(R(x_2))
        \end{equation}
        for all $x_1' \in \Post_{1, u_1}(R(x_2))$.
        By \defref{def:altsim}, there exists $u_2 \in \mathcal{U}_{2}$ such that
        $(R(x_2), x_2, u_1, u_2) \in R^e$, that is, for all $x_2' \in \Post_{2,u_2}(x_2)$,
        we have $R(x_2) \dyntou{u_1} R(x_2')$.
        By \eqref{eq:decrease} and \eqref{eq:lyap_proof}, we have
        \begin{align*}
            L_2(x_2) = L_1(R(x_2))
            & \ge L_1(R(x_2')) + c_1(R(x_2), u_1)\\
            & \ge L_2(x_2') + c_2(x_2, u_1).
        \end{align*}
    \end{myproof}
\end{theorem}

Receding horizon controllers may not reach the target due to their short-sighted nature.
This can be circumvented thanks to a Lyapunov function in several ways, two of which we recall in \propref{mpc1} and \propref{mpc2}.

\begin{algorithm}
  \SetAlgoLined
  \KwData{Initial state $x_0$, target set $\mathcal{X}_t$, horizon $\horizon{}$, Q-functions $(Q_k(x))_{k=0}^\infty$.}
  $k \gets 0$\;
  
  \While{$x_k \notin \mathcal{X}_t$}{
     $\tuples{u}{1}{\horizon} \in \argmin_{\tuple{u}{1}{\horizon} \in \mathcal{U}^\horizon} Q_k(x_k, \tuple{u}{1}{\horizon})$\;
     
     $u_{k+1} \gets u^\star_1$\;
     
     $x_{k+1} \in \Post_{u_{k+1}}(x_k)$\;
     
     $k \gets k + 1$\;
  }
  \Return{$u, x$}\;
  \caption{
    Receding horizon controller algorithm for a discrete-time control system as defined in \defref{def:discrete_time}.
  }
  \label{algo:mpc}
\end{algorithm}

The following proposition provides a classical condition for ensuring the convergence of
a model predictive controller \cite{mayne2001control}.

\begin{proposition}[\cite{mayne2001control}.]
  \label{prop:mpc1}
  Consider a discrete-time control system $S = (\mathcal{X}, \mathcal{U}, \dynto)$,
  a target set $\mathcal{X}_t \subseteq \mathcal{X}$,
  a set $\mathcal{X}_f \subseteq \mathcal{X}$ and
  a nonnegative cost function $c$.
  Let $L$ be a nonnegative Lyapunov function with cost function $c$ for $S$ in $\mathcal{X}_f$.
  Let $v_k = \min_{\tuple{u}{1}{\horizon} \in \mathcal{U}^\horizon} Q_k(x_k, \tuple{u}{1}{\horizon})$.
  If $Q_k = \QcV{L}$ for all $k \in \mathbb{N}$, and
  \begin{equation}
    \label{eq:delta}
    \delta = \inf_{x \in \mathcal{X} \setminus \mathcal{X}_t,\ u \in \mathcal{U}} c(x, u) > 0
  \end{equation}
  then \cref{algo:mpc} terminates in at most $v_0 / \delta$ iterations.
  \begin{myproof}
    We prove that after the $k$th iteration $v_k \le v_{k-1} - \delta$.
    Since $c$ and $L$ are nonnegative, $v_k$ is nonnegative for all $k$.
    Therefore, as $v_0$ is finite, this shows that \cref{algo:mpc} terminates in at most $v_0 / \delta$ iterations.

    Consider $u^\star$ found at the $k$th iteration of \cref{algo:mpc}.
    $Q_k(x_k, u^\star)$ is finite, hence $x' \in \mathcal{X}_f$ for all $x' \in \Post_{u^\star}(x_k)$.
    By \defref{def:lyap}, for all $x' \in \Post_{u^\star}(x_i)$, there exists $u'(x')$ such that $L(x'') \le L(x') - c(x', u')$ for all $x'' \in \Post_{u'}(x')$.
    Therefore, for all $x_{k+1} \in \Post_{u^\star_1}(x_k)$,
    $Q_{k+1}(x_{k+1}, \tuplecat{\tuples{u}{2}{\horizon}}{u'}) \le Q_{k+1}(x_{k+1}, \tuples{u}{2}{\horizon{}})$.
    And finally, by \eqref{eq:delta}, we have $Q_{k+1}(x_{k+1}, \tuples{u}{2}{\horizon{}}) \le Q_k(x_k, u^\star) - \delta$.
  \end{myproof}
\end{proposition}


The following proposition generalizes \cite[Theorem~5.4]{gol2015temporal} where the Lyapunov function is called ``distance function''.
This distance function is computed from the Lyapunov function of an alternatingly simulated system
that is constructed with \cite[Algorithm~2]{gol2013language}.

\begin{proposition}
  \label{prop:mpc2}
  Consider a deterministic discrete-time control system $S = (\mathcal{X}, \mathcal{U}, \dynto)$,
  a target set $\mathcal{X}_t \subseteq \mathcal{X}$,
  a set $\mathcal{X}_f \subseteq \mathcal{X}$
  and the cost function $c$ such that for all $x \in \mathcal{X}$ and $u \in \mathcal{U}$,
  $c(x, u) = 0$ if $x \in \mathcal{X}_t$ and $c(x, u) = 1$ otherwise.
  Let $L$ be a Lyapunov function with cost function $c$ for $S$ in $\mathcal{X}_f$
  such that $L(x) = 0$ if $x \in \mathcal{X}_t$.
  Suppose there is $\deadline \in \mathbb{N}$ and value functions $\tuple{V}{0}{\deadline}$ such that
  $V_i(x)$ is finite if and only if $i \ge L(x)$ and
  $Q_k(x, \tuple{u}{1}{\horizon}) = \QcV{V_{\max(0, \deadline - k - \horizon)}}(x, \tuple{u}{1}{\min(\horizon, \deadline - k)})$
  for $k = 1, 2, \ldots, \deadline$, $x \in \mathcal{X}$ and $\tuple{u}{1}{\horizon} \in \mathcal{U}^\horizon$.
  Let $v_k = \min_{\tuple{u}{1}{\horizon} \in \mathcal{U}^\horizon} Q_k(x_k, \tuple{u}{1}{\horizon})$.
  If $v_0$ is finite, then
  then \cref{algo:mpc} terminates in at most $\deadline$ iterations.
  \begin{myproof}
    We prove by induction that $v_k$ is finite
    for $k = 0, \ldots, \deadline$.
    In particular, the finiteness of $v_{\deadline}$ implies the statement of the proposition.

    If $v_k$ is finite, then there exists $\tuples{u}{1}{\horizon{}}$ such that
    $Q_k(x_k, \tuples{u}{1}{\horizon{}})$ is finite, hence $L(x') \le \max(0, \deadline - k - \horizon)$ for the unique $x' \in \Post_{\tuples{u}{1}{\horizon{}}}(x_k)$.
    By \defref{def:lyap}, there exists $u'$ such that $L(x'') \le \max(0, \deadline - k - \horizon - 1)$ for the unique $x'' \in \Post_{u'}(x')$.
    Therefore, for all $x_{k+1} \in \Post_{u^\star_1}(x_k)$, $Q_{k+1}(x_{k+1}, \tuplecat{\tuples{u}{2}{\horizon{}}}{u'})$ is finite, and so is $v_{k+1}$.
  \end{myproof}
\end{proposition}



\section{\BnB{} algorithm}
\label{sec:mip}

In this section, we show how the concepts of Lyapunov functions and Bellman-like Q-functions can be exploited by a \bnb{} algorithm.
In \cref{sec:Qlearning}, we show that the Bellman-like Q-function can be further refined during the optimization by learning it only along the optimal trajectory as approximating it over the whole state-space is not scalable.

For this section we use the following definition of hybrid systems.
As it is a special case of \defref{def:discrete_time}, it allows to reuse the results of the previous section.

\begin{definition}[\Hsyst{}]\label{def:controlled_hybrid_automaton}
    A \emph{\hsyst{}} is defined as a triple
    $\Hsys = (\mathcal{Q} \times \mathcal{X}, \mathcal{V} \times \mathcal{U}, \dynto)$
    where $\mathcal{Q}$ is the finite set representing the discrete state-space,
    $\mathcal{X} \subseteq \mathbb{R}^{n_x}$ represents the continuous state-space,
    $\mathcal{V}$ is the finite set of discrete control inputs,
    $\mathcal{U} \subseteq \mathbb{R}^{n_u}$ is the set of continuous control inputs
    and $\dynto$ is the subset of
    transitions $((q, x), (v, u), (q', x'))$ such that the system can reach $q' \in \mathcal{Q}, x' \in \mathcal{X}$ from $q \in \mathcal{Q}, x \in \mathcal{X}$ with inputs $v \in \mathcal{V}, u \in \mathcal{U}$.
    We denote $((q, x), (v, u), (q', x')) \in \dynto$ as $(q, x) \dyntou{v, u} (q', x')$. 
\end{definition}


The optimal control problem is formally defined as follows.

\begin{problem}
   \label{prob:optimal_control}
   Consider a \hsyst{} $\Hsys$ as defined in \defref{def:controlled_hybrid_automaton}.
   The optimal control problem for $\Hsys$ with initial \states{} $q_0 \in \mathcal{Q}, x_0 \in \mathcal{X}$,
   target set $\mathcal{X}_t$
   and cost function
   $c$
   is defined as
   the optimization problem:
   \begin{equation}
       \label{eq:optimal_control}
       \inf_{
           l \in \mathbb{N},\ \tuple{v}{1}{l} \in \mathcal{V}^l,\ \tuple{u}{1}{l} \in \mathcal{U}^l
       } \Qvu[\QcV{V_0}]{q_0}{x_0}{\tuple{v}{1}{l}}{\tuple{u}{1}{l}}
   \end{equation}
   where $V_0(q, x) = 0$ for all $(q, x) \in \mathcal{X}_t$ and $V_0(q, x) = \infty$ otherwise.
\end{problem}

Given a Bellman-like Q-function $Q$, we define the
Q-function $\Qfunv$ that is only parametrized by the discrete input $\tuple{v}{1}{k}$ as:
\begin{equation}
   \label{eq:Bellman-likeQv}
   \Qv{q}{x}{\tuple{v}{1}{k}} = \min_{\tuple{u}{1}{k} \in \mathcal{U}^k} \Qvu{q}{x}{\tuple{v}{1}{k}}{\tuple{u}{1}{k}}.
\end{equation}
The following proposition shows that this is a Bellman-like Q-function as well.

\begin{proposition}
   \label{prop:Bellman-likeQv}
   Consider a \hsyst{} $\Hsys$ as defined in \defref{def:controlled_hybrid_automaton}.
   If $Q$ is a Bellman-like Q-function of $\Hsys$ with cost function $c$,
   then the Q-function $\Qfunv$ defined by \eqref{eq:Bellman-likeQv}
   is a Bellman-like Q-function of $\Hsys$ with cost function $c$.
   \begin{myproof}
      For given states $q, x$ and discrete inputs $\tuple{v}{1}{l}$, for any $k < l$
      and $\tuple{u}{1}{l} \in \mathcal{U}^l$, by \defref{def:Bellman-likeQ}, we have
      $\Qvu{q}{x}{\tuple{v}{1}{k}}{\tuple{u}{1}{k}} \le
      \Qvu{q}{x}{\tuple{v}{1}{l}}{\tuple{u}{1}{l}}$.
      Combining this inequality with \eqref{eq:Bellman-likeQv}, we obtain
      \[
          \Qv{q}{x}{\tuple{v}{1}{l}}
          = \min_{\tuple{u}{1}{l} \in \mathcal{U}^l} \Qvu{q}{x}{\tuple{v}{1}{l}}{\tuple{u}{1}{l}}
          \ge \min_{\tuple{u}{1}{k} \in \mathcal{U}^k} \Qvu{q}{x}{\tuple{v}{1}{k}}{\tuple{u}{1}{k}}
          = \Qv{q}{x}{\tuple{v}{1}{k}}.
      \]
   \end{myproof}
\end{proposition}

As shown in \cite{bemporad1999control}, \cref{prob:optimal_control} can be formulated as a Mixed Integer Quadratic Program (MIQP) and then solved by generic MIQP solvers.
On the other hand, we show in the remaining of this section that \cref{algo:BnB} can incorporate the information gathered in the computation of a Lyapunov and Bellman-like value functions, as well as refine these functions during the optimization.
The aim is to enable the \bnb{} algorithm
to better exploit the structure of the problem than a generic MIQP solver.

\begin{algorithm}[htbp]
    \SetAlgoLined
    \KwData{
    Initial \states{} $q_0 \in \mathcal{Q}$, $x_0 \in \mathcal{X}$, a heuristic function $h$, a target set $\mathcal{X}_t$,
    a Bellman-like Q-function $Q$, an upper bound function $\beta$ and $\lenub \in \mathbb{N}$ satisfying \cref{assu:deadline}.}
    \PrintSemicolon
    $\bar{\beta} \gets \infty$\;
    
    $\mathcal{N} \gets \{\emptytuple{}\}$\;

    \While{$\mathcal{N} \neq \varnothing$}{
        $\tuple{v}{1}{l} \gets h(\mathcal{N})$\;
        
        $\mathcal{N} \gets \mathcal{N} \setminus \{\tuple{v}{1}{l}\}$\;

        \If{$\Qv{q_0}{x_0}{\tuple{v}{1}{l}} \le \bar{\beta}$ and $l < \lenub$}{
            \For{$v' \in \mathcal{V}$}{
                $\tuple{\hat{v}}{1}{\hat{l}}, \tuple{\hat{u}}{1}{\hat{l}}, \hat{\beta} \gets \beta(q_0, x_0, \tuplecat{\tuple{v}{1}{l}}{v'})$\;
                
                \If{$\hat{\beta} < \bar{\beta}$}{
                   $\tuple{\bar{v}}{1}{\bar{l}}, \tuple{\bar{u}}{1}{\bar{l}}, \bar{\beta} \gets \tuple{\hat{v}}{1}{\hat{l}}, \tuple{\hat{u}}{1}{\hat{l}}, \hat{\beta}$\;
                   
                   
                }
                $\mathcal{N} \gets \mathcal{N} \cup \{\tuplecat{\tuple{v}{1}{l}}{v'}\}$\;
            }
        }
    }
    \Return{$\tuple{\bar{v}}{1}{\bar{l}}, \tuple{\bar{u}}{1}{\bar{l}}$}\;
    \caption{
      \BnB{} algorithm for \cref{prob:optimal_control}.
      The set $\mathcal{N}$ represents the set of nodes of the search tree for which
      subtrees still need to be explored.
      The heuristic function determines which node is considered next,
      two different heuristics are discussed in \cref{sec:example}.
    }
    \label{algo:BnB}
\end{algorithm}

\subsection{Convergence and optimality}
\label{sec:opt}

In this section is discussed the convergence and optimality of the algorithm.
For optimality, we need to ensure that the condition ``$\Qv{q_0}{x_0}{\tuple{v}{1}{l}} \le \bar{\beta}$''
in the algorithm does not exclude any optimal solution.
To this end, we start by proving in \lemref{lower} that $\Qv{q_0}{x_0}{\tuple{v}{1}{l}}$ gives a lower bound to \eqref{eq:optimal_control},
provided that $\Qfunv$ is a Bellman-like Q-function.


\begin{lemma}
  \label{lem:lower}
  Consider \cref{prob:optimal_control} for a \hsyst{} $\Hsys$
  with initial \states{} $q_0, x_0$, cost function $c$, and a Bellman-like Q-function
  $Q$ for $\Hsys$ with cost function $c$.
  If
  $\Qvu{q_0}{x_0}{\tuple{v}{1}{l}}{\tuple{u}{1}{l}}
  \le \Qvu[\QcV{V_0}]{q_0}{x_0}{\tuple{v}{1}{l}}{\tuple{u}{1}{l}}$
  for any $\tuple{v}{1}{l} \in \mathcal{V}^l$ and $\tuple{u}{1}{l} \in \mathcal{U}^l$,
  then $\Qv{q_0}{x_0}{\tuple{v}{1}{k}}$ is a lower bound to \eqref{eq:optimal_control} for any $\tuple{v}{1}{k} \in \mathcal{V}^k$.
  \begin{myproof}
     Consider $k, l \in \mathbb{N}$ such that $k \le l$ and $\tuple{v}{1}{l} \in \mathcal{V}^l$.
     By \propref{Bellman-likeQv}, $\Qv{q_0}{x_0}{\tuple{v}{1}{k}}$ is a Bellman-like Q-function hence
     \defref{def:Bellman-like} ensures that $\Qv{q_0}{x_0}{\tuple{v}{1}{k}} \le \Qv{q_0}{x_0}{\tuple{v}{1}{l}}$.
     By \eqref{eq:Bellman-likeQv}, we have $\Qv{q_0}{x_0}{\tuple{v}{1}{l}} \le \Qvu{q_0}{x_0}{\tuple{v}{1}{l}}{\tuple{u}{1}{l}}$.
     The statement of the lemma ensures that
     $\Qv{q_0}{x_0}{\tuple{v}{1}{k}} \le \Qvu[\QcV{V_0}]{q_0}{x_0}{\tuple{v}{1}{l}}{\tuple{u}{1}{l}}$
     which concludes the proof.
  \end{myproof}
\end{lemma}

To ensure the convergence of the algorithm, the following assumption excludes pathological
optimal control problems that only admit arbitrarily long optimal solutions.

\begin{assumption}
  \label{assu:deadline}
  There exists $\lenub \in \mathbb{N}$ and an optimal solution $\tuple{v}{1}{\lenub}, \tuple{u}{1}{\lenub}$ of \cref{prob:optimal_control}. 
\end{assumption}

\begin{remark}
  \label{rem:deadline}
  \cref{assu:deadline} ensures that there exists an optimal solution of finite length.
  A more conservative alternative to \cref{assu:deadline} would be to assume that the cost function $c$ is lower bounded by
  a positive number and that there is an upper bound to the optimal cost of the optimal control problem.
\end{remark}

The following result ensures both the convergence and optimality of \cref{algo:BnB}.

\begin{theorem}
  \label{theo:BnB}
  Consider \cref{prob:optimal_control} for a deterministic \hsyst{} $\Hsys$
  with initial \states{} $q_0, x_0$, cost function $c$, a Bellman-like Q-function
  $Q$ for $\Hsys$ with cost function $c$, and an upper bound function $\beta$.
  Assume that
  $\Qvu{q_0}{x_0}{\tuple{v}{1}{l}}{\tuple{u}{1}{l}}
  \le \Qvu[\QcV{V_0}]{q_0}{x_0}{\tuple{v}{1}{l}}{\tuple{u}{1}{l}}$
  for any $\tuple{v}{1}{l} \in \mathcal{V}^l$ and $\tuple{u}{1}{l} \in \mathcal{U}^l$
  and that $\beta$ either returns $\emptytuple{}, \emptytuple{}, \infty$ or
  $v, u, \Qvu[\QcV{V_0}]{q_0}{x_0}{\tuple{v}{1}{l}}{\tuple{u}{1}{l}}$
  with finite $\Qvu[\QcV{V_0}]{q_0}{x_0}{\tuple{v}{1}{l}}{\tuple{u}{1}{l}}$.
  Then \cref{algo:BnB} returns an optimal solution of \cref{prob:optimal_control}.
  \begin{myproof}
    As $\mathcal{V}$ is finite, the condition ``$l < \lenub$'' of \cref{algo:BnB}
    guarantees that the algorithm terminates.
    To prove that \cref{algo:BnB} returns an optimal solution, we prove that \cref{algo:BnB} calls
    $\beta(q_0, x_0, \tuplecat{\tuple{v}{1}{l}}{v'})$ for $\tuplecat{\tuple{v}{1}{l}}{v'}$ equal to the discrete transitions
    of one of the optimal trajectories whose existence is guaranteed by \cref{assu:deadline}.
    Indeed, for any prefix $\tuple{v}{1}{k}$ of $\tuple{v}{1}{l}$, the condition ``$k < \lenub$'' is satisfied as
    the length of $\tuplecat{\tuple{v}{1}{l}}{v'}$ is at most $\lenub$
    and the condition ``$\Qv{q_0}{x_0}{\tuple{v}{1}{l}} \le \bar{\beta}$'' is satisfied by \lemref{lower}.
  \end{myproof}
\end{theorem}

\subsection{Obtaining upper bounds}

\Cref{algo:BnB} is parametrized by a function $\beta$ responsible
to provide upper bounds.
As discussed in \cref{sec:opt}, the only necessary condition on $\beta$ for the convergence and optimality of
\cref{algo:BnB} is that $\beta$ should either return nothing or a feasible solution of \cref{prob:optimal_control}.
However, a good algorithm for $\beta$, i.e. that provides a feasible solution of low cost,
can have a dramatic impact on the efficiency of \cref{algo:BnB} as it allows it to prune significant parts of the search tree.

In this section, we introduce a candidate for $\beta$ as \cref{algo:upper}.
\Cref{algo:upper} searches for feasible trajectories given a fixed prefix of discrete inputs.
Its ability to return feasible solutions highly depends on the size of the target set $\mathcal{X}_f$ and the horizon $\horizon{}$.
However, computing a Lyapunov function with a larger target set $\mathcal{X}_f$ requires more offline computation,
while a larger horizon $\horizon{}$ requires more online computation in \cref{algo:upper}.
This increase in computational effort might result in better pruning for \cref{algo:BnB},
hence there is a compromise to reach between a computationally cheap $\beta$ function that often provides
a costly feasible solution or no feasible solution at all, and
a computationally expensive $\beta$ function that will quickly find good feasible solutions by pruning large parts of the search tree.


\begin{algorithm}
  \SetAlgoLined
  \KwData{A deterministic \hsyst{} $\Hsys$, initial \states{} $q_0 \in \mathcal{Q}, x_0 \in \mathcal{X}$,
  discrete input $\tuple{v}{1}{k}$,
  horizon $\horizon{}$, deadline $\deadline{}$ and
  Q-functions $(Q_k)_{k=0}^{\deadline{}}$.}
  $Z(q, x) = 0$\;

  $\tuple{u}{1}{k} \gets \argmin_{\tuple{u}{1}{k} \in \mathcal{U}^k} \Qvu[\QcV{Z}]{q_0}{x_0}{\tuple{v}{1}{k}}{\tuple{u}{1}{k}}$\;
  
  \If{$\Qvu[\QcV{Z}]{q_0}{x_0}{\tuple{v}{1}{k}}{\tuple{u}{1}{k}} = \infty$}{
    \Return{$\emptytuple{}, \emptytuple{}, \infty$}\;
  }
  Let $(q_k, x_k)$ be the unique pair such that $(q_0, x_0) \dyntou{\tuple{v}{1}{k}, \tuple{u}{1}{k}} (q_k, x_k)$\;
  
  \eIf{$\min_{\tuple{v}{1}{\horizon} \in \mathcal{V}^\horizon, \tuple{u}{1}{\horizon} \in \mathcal{U}^\horizon}
  \Qvu[Q_k]{q_k}{x_k}{\tuple{v}{1}{\horizon}}{\tuple{u}{1}{\horizon}}$ is finite}{
    $(\tuplep{v}{1}{l'}, \tuplep{u}{1}{l'}), (q', x') \gets {}$ Algorithm~1 with $(q_k, x_k)$, $\horizon{}$, $\deadline{} - k$, $(Q_i)_{i=k+1}^{\deadline{}}$\;
    
    \Return{$\tuplecat{\tuple{v}{1}{k}}{\tuplep{v}{1}{l'}}, \tuplecat{\tuple{u}{1}{k}}{\tuplep{u}{1}{l'}},
      \Qvu[\QcV{V_0}]{q_0}{x_0}{\tuplecat{\tuple{v}{1}{k}}{\tuplep{v}{1}{l'}}}{\tuplecat{\tuple{u}{1}{k}}{\tuplep{u}{1}{l'}}}$}\;
  }{
    \Return{$\emptytuple{}, \emptytuple{}, \infty$}\;
  }
  \caption{
  Upper bound algorithm that can be used as $\beta$ function for \cref{algo:BnB}.
  }
  \label{algo:upper}
\end{algorithm}

\begin{proposition}
   \label{prop:upper}
   Consider \cref{prob:optimal_control} for a deterministic \hsyst{} $\Hsys$
   with initial \states{} $q_0, x_0$ and cost function $c$.
   If the Q-functions $(Q_k(x))_{k=0}^{\deadline{}}$ satisfy either the assumptions of \propref{mpc1} or \propref{mpc2},
   then \cref{algo:upper} either returns $\emptytuple{}, \emptytuple{}, \infty$ or
   $\tuple{v}{1}{l}, \tuple{u}{1}{l}, \Qvu[\QcV{V_0}]{q_0}{x_0}{\tuple{v}{1}{l}}{\tuple{u}{1}{l}}$
   with finite $\Qvu[\QcV{V_0}]{q_0}{x_0}{\tuple{v}{1}{l}}{\tuple{u}{1}{l}}$.
   \begin{myproof}
      Suppose that \cref{algo:upper} returns
      $\tuple{v}{1}{l}, \tuple{u}{1}{l}, \Qvu[\QcV{V_0}]{q_0}{x_0}{\tuple{v}{1}{l}}{\tuple{u}{1}{l}}$
   with finite $\Qvu[\QcV{V_0}]{q_0}{x_0}{\tuple{v}{1}{l}}{\tuple{u}{1}{l}}$.
      By \propref{mpc1} and \propref{mpc2}, $(q_k, x_k) \dyntou{\tuplep{v}{1}{l'}, \tuplep{u}{1}{l'}} (q', x')$
      with $(q', x') \in \mathcal{X}_t$.
      Therefore, as $\Qvu[\QcV{Z}]{q_0}{x_0}{\tuple{v}{1}{k}}{\tuple{u}{1}{k}}$ is finite, $(q_0, x_0) \dyntou{\tuple{v}{1}{k}, \tuple{u}{1}{k}} (q_k, x_k)$,
      hence $(q_0, x_0) \dyntou{\tuplecat{\tuple{v}{1}{k}}{\tuplep{v}{1}{l'}}, \tuplecat{\tuple{u}{1}{k}}{\tuplep{u}{1}{l'}}} (q', x')$.
   \end{myproof}
\end{proposition}

\subsection{Q-learning}
\label{sec:Qlearning}

The difference between the value of the Q-function provided to \cref{algo:BnB} with the actual minimal cost of \cref{prob:optimal_control}
has a significant impact on the performance of \cref{algo:BnB}.
As discussed in \cite{bouchatreinforcement, menta2020learning}, the function providing this minimal cost
is in general nonlinear and nonconvex.
As a matter of fact, a Q-function with a small such difference over the whole state space may not be computable in a reasonable amount of time.

To circumvent this issue, \cite{bouchatreinforcement} suggests a reinforcement learning approach
to generate an approximation of the Q-function that is close to the actual minimal cost near the optimal trajectory
of the problem.
As the optimal trajectory is unknown, the approach employed by \cite{bouchatreinforcement},
which is classical in Stochastic Programming~\cite{birge2011introduction},
consists in alternating between a \emph{forward pass} and a \emph{backward pass}.
The forward pass computes a feasible trajectory with a receding horizon controller
using the current Q-function as terminal cost.
Starting from the end of the trajectory, the backward pass generates new cuts for the Q-function corresponding to
each transition using the state at each step of the trajectory.

This backward pass can be used to refine the Q-functions along feasible trajectories found by $\beta$ in
\cref{algo:BnB}.
Given such trajectory $(\hat{q}_0, \hat{x}_0) \dyntou{\hat{v}_1, \hat{u}_1} (\hat{q}_1, \hat{x}_1) \dyntou{\hat{v}_2, \hat{x}_2} \cdots \dyntou{\hat{v}_l, \hat{u}_l} (\hat{q}_l, \hat{x}_l)$
and a Q-function $\Qvu{q}{x}{\tuple{v}{1}{l}}{\tuple{u}{1}{l}}$ defined for $l \ge k$ with \eqref{eq:Qrec}, this learning phase consists in
replacing $Q$ by $\max(Q, Q_j')$ for $j = l-k, l-k-1, \ldots, 1$ where $Q_j'$ is computed as follows.
The affine lower approximation $Q'_j$ of the function
$\mathcal{T}_c \Qvu{q}{x}{\tuple{v}{1}{k}}{\tuple{u}{1}{k}}
$ is obtained using a feasible solution of the dual of the problem
$\min_{v',u'} \Qvu{\hat{q}_j}{\hat{x}_j}{\tuplecat{\tuple{\hat{v}}{j+1}{j+k}}{v'}}{\tuplecat{\tuple{\hat{u}}{j+1}{j+k}}{u'}}$.
As shown in the following proposition, the set of Bellman-like Q-functions is invariant under this operation.

\begin{proposition}
   Consider \cref{prob:optimal_control} for a deterministic \hsyst{} $\Hsys$
   with cost function $c$
   and a Bellman-like Q-function $Q$ for $\Hsys$ with $c$.
   If $\Qvu[Q']{q}{x}{\tuple{v}{1}{k}}{\tuple{u}{1}{k}} \le \mathcal{T}_c \Qvu{q}{x}{\tuple{v}{1}{k}}{\tuple{u}{1}{k}}$, then 
   $\max(Q, Q')$ is a Bellman-like Q-function.
   \begin{myproof}
      We prove this proposition by contradiction.
      Consider $q, x, \tuple{v}{1}{s+k+1}, \tuple{u}{1}{s+k+1}$ such that
      $$\Qvu[\max(Q, Q')]{q}{x}{\tuple{v}{1}{s+k}}{\tuple{u}{1}{s+k}}
      > \Qvu[\max(Q, Q')]{q}{x}{\tuple{v}{1}{s+k+1}}{\tuple{u}{1}{s+k+1}}.$$
      Let $q', x'$ be such that $(q, x) \dyntou{\tuple{v}{1}{s}, \tuple{u}{1}{s}} (q', x')$, $\tuplep{v}{1}{k+1} = \tuple{v}{s+1}{s+k+1}$, and $\tuplep{u}{1}{k+1} = \tuple{u}{s+1}{s+k+1}$.
      We have
      $ \Qvu[\max(Q, Q')]{q'}{x'}{\tuplep{v}{1}{k  }}{\tuplep{u}{1}{k  }}
      > \Qvu[\max(Q, Q')]{q'}{x'}{\tuplep{v}{1}{k+1}}{\tuplep{u}{1}{k+1}}$
      which implies in particular that
      $ \Qvu[\max(Q, Q')]{q'}{x'}{\tuplep{v}{1}{k  }}{\tuplep{u}{1}{k  }}
      >              \Qvu{q'}{x'}{\tuplep{v}{1}{k+1}}{\tuplep{u}{1}{k+1}}.$
      As $Q$ is a Bellman-like Q-function,
      $\Qvu{q'}{x'}{\tuplep{v}{1}{k  }}{\tuplep{u}{1}{k  }} \le
       \Qvu{q'}{x'}{\tuplep{v}{1}{k+1}}{\tuplep{u}{1}{k+1}}$
      hence we have
      $\Qvu[Q']{q'}{x'}{\tuplep{v}{1}{k  }}{\tuplep{u}{1}{k  }} >
           \Qvu{q'}{x'}{\tuplep{v}{1}{k+1}}{\tuplep{u}{1}{k+1}}$
      which is in contradiction with the statement of the proposition.
   \end{myproof}
\end{proposition}

\section{Numerical example}
\label{sec:example}

In this section, we illustrate the algorithms developed in this paper on the double integrator dynamics example introduced in \cite[Example~VIII.A]{gol2013language}.
The deterministic discrete-time hybrid control automaton $(\mathcal{Q} \times [-10, 1.85]^2, \mathcal{V} \times [-2, 2], \dynto)$ is such that a transition $(q, x) \dyntou{v, u} (q', x')$
occurs if
\[
  x' = \begin{bmatrix}1 & 1\\0 & 1\end{bmatrix} x + \begin{bmatrix}0.5\\1.0\end{bmatrix},
\]
in addition to logical constraints ensuring that a feasible trajectory goes through either square
$\mathcal{A}$ or $\mathcal{B}$ before reaching the target square $\mathcal{X}_t$. The squares are represented in \Cref{fig:plots}.
See \cite[Example~VIII.A]{gol2013language} for more details on the definition of the system.
Finally, let the cost $c((q, x), (v, u)) = u_1^2$ if $q = \modeT$, which is equivalent to $x \in \mathcal{X}_t$, and $c((q, x), (v, u)) = u_1^2 + 1$ otherwise.

We benchmark the number of iterations of the branch and bound algorithm with no Lyapunov function
and horizon $\horizon = 0$, i.e. upper bounds are only obtained when the candidate $\tuple{v}{1}{l}$ of \cref{algo:BnB}
has $\modeT$ as final discrete state.
The heuristic $h$ is a depth-first heuristic, i.e., it selects the candidate $\tuple{v}{1}{l}$ with the largest $k$
and breaks ties by selecting the one with the smallest lower bound $\Qv{q_0}{x_0}{\tuple{v}{1}{l}}$.
Note that, as we have no Lyapunov function and a horizon $\horizon=0$, a breadth-first heuristic
would not be able to prune much of the search tree as $\bar{\beta}$ would remain infinite for
most of the iterations.

In order to study the generalization of the Q-function
we consider two instances of the optimal control problem: $\mathcal{I}_1$ with $x_0 = [1.5, -2.5], \lenub = 9$ and $\mathcal{I}_2$ with $x_0 = [1, -6], \lenub = 11$.

We analyze the behavior of \cref{algo:BnB} with two different Bellman-like Q-functions: $Q_1$ and $Q_2$.
The first one
is the trivial $Q_1(x, u) = -\infty$ that corresponds to no lower bound
hence no pruning in the branch and bound algorithm.
The second one, $Q_2$, is obtained by applying \cref{theo:bellman}
on the alternating simulation $(\mathcal{Q}, \mathcal{V}, \dynto')$ such that
$q \dyntou{v}' q'$ if there exists $x, x' \in [-10, 1.85]^2$, $u \in [-2, 2]$ such that
$(q, x) \dyntou{v, u} (q', x')$.
The number of iterations for different choices of Bellman-like Q-function is given in \cref{tab:learn}.
Feasible trajectories found by $\beta$ are given in \cref{fig:plots}.

\begin{table}[htbp]
    \centering
    \begin{tabular}{c|c|c}
        Bellman-like Q-function & Number of iterations for $\mathcal{I}_1$ & Number of iterations for $\mathcal{I}_2$\\
        \hline
        $Q_1$ & 197234 & 9388410\\
        $Q_2$ & 1111 & 85\\
        $Q_2$ with learning & 871 & 85\\
        $Q_2$ learning on $\mathcal{I}_1$ & 761 & 85\\
        $Q_2$ learning on $\mathcal{I}_2$ & 880 & 85\\
        $Q_2$ learning on $\mathcal{I}_1$ and $\mathcal{I}_2$ & 747 & 85\\
    \end{tabular}
    \caption{
      Number of iterations of \cref{algo:BnB} with input parameters described in \cref{sec:example}
      for the instances $\mathcal{I}_1$ and $\mathcal{I}_2$ with different Bellman-like Q-functions.
      We observe that the number of iterations for $\mathcal{I}_2$ is drastically reduced%
      . When using the abstraction-based Q-function $Q_2$, it is divided by more than $10^5$.
      It is also significantly reduced for $\mathcal{I}_1$, by more than $10^3$, and it is further
      decreased thanks to the learning approach detailed in \cref{sec:Qlearning}.
      We also note that the Bellman-like Q-function learned on $\mathcal{I}_2$
      generalize well for $\mathcal{I}_1$, with only 880 iterations, and
      the Q-function learned on the same instance prunes even more nodes, 
      as there is only 761 iterations.
      Combining the cuts learned from both $\mathcal{I}_1$ and $\mathcal{I}_2$ prunes even more nodes with only 747 iterations needed.
      This shows the ability of the Q-learning approaches to both generalize and combine knowledge about the Bellman-like Q-functions.
    }
    \label{tab:learn}
\end{table}


\section{Conclusion}


The size of abstractions that can simulate the behavior of a hybrid system
with enough accuracy
in the whole state-space typically grows exponentially with the dimension of the systems.
However, as we show in \cref{sec:abstraction},
any abstraction can provide a Lyapunov function over some set $\mathcal{X}_f$ or a Bellman-like value function.
Of course, the coarser the abstraction, the smaller the set $\mathcal{X}_f$, and the larger the gap between the
value of the Q-function provided to \cref{algo:BnB} and the actual minimal cost of \cref{prob:optimal_control}.
Nevertheless, we show in \cref{sec:mip} that this information can be leveraged by a branch and bound algorithm.
Moreover, as illustrated by our numerical example in \cref{sec:example}, even the Bellman-like value function obtained from a rather coarse abstraction allows drastic pruning of the search tree of the \bnb{} algorithm.

While the computation of a global Lyapunov function or a good approximation of the minimal cost of \cref{prob:optimal_control} 
is not tractable for hybrid systems in general,
we can still aim at computing a local Lyapunov function and a Bellman-like value function that are sufficient
for an effective pruning of the search tree.
For this purpose, it seems appropriate to guide the improvements of these functions using feasible trajectories
found during the algorithm.
Often, this will in practice enhance the refinement of these functions in the appropriate regions of the state-space.
As we show in \cref{sec:example}, the refinement of the Bellman-like value function obtained
after \cref{algo:BnB} can be reused for solving similar optimal control problem.

Several key research directions of this approach are left as future work.
This includes a detailed complexity analysis of the algorithm with, in particular, the complexity of computing the cut in \cref{sec:Qlearning}.
A second line of research is the iterative refinement of the abstractions used to compute the Lyapunov
and Bellman-like value functions throughout the algorithms.



\acks{RJ is a FNRS honorary Research Associate. This project has received funding from the European Research Council (ERC) under the European Union's Horizon 2020 research and innovation programme under grant agreement No 864017 - L2C. RJ is also supported by the Walloon Region, the Innoviris Foundation, and the FNRS (Chist-Era Druid-net).}

\bibliography{biblio}

\begin{thebibliography}{14}
\providecommand{\natexlab}[1]{#1}
\providecommand{\url}[1]{\texttt{#1}}
\expandafter\ifx\csname urlstyle\endcsname\relax
  \providecommand{\doi}[1]{doi: #1}\else
  \providecommand{\doi}{doi: \begingroup \urlstyle{rm}\Url}\fi

\bibitem[Bemporad and Morari(1999)]{bemporad1999control}
Alberto Bemporad and Manfred Morari.
\newblock Control of systems integrating logic, dynamics, and constraints.
\newblock \emph{Automatica}, 35\penalty0 (3):\penalty0 407--427, 1999.

\bibitem[Bemporad et~al.(2002)Bemporad, Morari, Dua, and
  Pistikopoulos]{bemporad2002explicit}
Alberto Bemporad, Manfred Morari, Vivek Dua, and Efstratios~N Pistikopoulos.
\newblock The explicit linear quadratic regulator for constrained systems.
\newblock \emph{Automatica}, 38\penalty0 (1):\penalty0 3--20, 2002.

\bibitem[Birge and Louveaux(2011)]{birge2011introduction}
John~R Birge and Francois Louveaux.
\newblock \emph{{I}ntroduction to stochastic programming}.
\newblock Springer Science \& Business Media, 2011.

\bibitem[Blondel and Tsitsiklis(2000)]{blondel2000survey}
Vincent~D Blondel and John~N Tsitsiklis.
\newblock {A} survey of computational complexity results in systems and
  control.
\newblock \emph{Automatica}, 36\penalty0 (9):\penalty0 1249--1274, 2000.

\bibitem[Bouchat and Jungers(2020)]{bouchatreinforcement}
Jean Bouchat and Rapha{\"e}l~M Jungers.
\newblock Reinforcement learning for the optimal control of hybrid systems.
\newblock Master's thesis, UCLouvain, 2020.

\bibitem[Gol et~al.(2013)Gol, Lazar, and Belta]{gol2013language}
Ebru~Aydin Gol, Mircea Lazar, and Calin Belta.
\newblock Language-guided controller synthesis for linear systems.
\newblock \emph{IEEE Transactions on Automatic Control}, 59\penalty0
  (5):\penalty0 1163--1176, 2013.

\bibitem[Gol et~al.(2015)Gol, Lazar, and Belta]{gol2015temporal}
Ebru~Aydin Gol, Mircea Lazar, and Calin Belta.
\newblock Temporal logic model predictive control.
\newblock \emph{Automatica}, 56:\penalty0 78--85, 2015.

\bibitem[Legat et~al.(2020{\natexlab{a}})Legat, Bouchat, and
  Jungers]{legat2020abstractionCO}
Beno{\^\i}t Legat, Jean Bouchat, and Raphaël~M. Jungers.
\newblock {Abstraction-based branch and bound approach to Q-learning for hybrid
  optimal control}.
\newblock \url{https://www.codeocean.com/}, November 2020{\natexlab{a}}.

\bibitem[Legat et~al.(2020{\natexlab{b}})Legat, Bouchat, and
  Jungers]{legat2020abstractionarXiv}
Beno{\^\i}t Legat, Jean Bouchat, and Raphaël~M. Jungers.
\newblock {Abstraction-based branch and bound approach to Q-learning for hybrid
  optimal control}.
\newblock \emph{ArXiv e-prints}, November 2020{\natexlab{b}}.

\bibitem[Legat et~al.(2020{\natexlab{c}})Legat, Dowson, Garcia, and
  Lubin]{legat2020mathoptinterface}
Beno\^{\i}t Legat, Oscar Dowson, Joaquim~Dias Garcia, and Miles Lubin.
\newblock Mathoptinterface: a data structure for mathematical optimization
  problems.
\newblock \emph{arXiv preprint arXiv:2002.03447}, 2020{\natexlab{c}}.

\bibitem[Mayne(2001)]{mayne2001control}
David~Q Mayne.
\newblock Control of constrained dynamic systems.
\newblock \emph{European Journal of Control}, 7\penalty0 (2-3):\penalty0
  87--99, 2001.

\bibitem[Menta et~al.(2020)Menta, Warrington, Lygeros, and
  Morari]{menta2020learning}
Sandeep Menta, Joseph Warrington, John Lygeros, and Manfred Morari.
\newblock Learning solutions to hybrid control problems using benders cuts.
\newblock In \emph{Learning for Dynamics and Control (L4DC) 2020}, 2020.

\bibitem[Stellato et~al.(2020)Stellato, Banjac, Goulart, Bemporad, and
  Boyd]{osqp}
B.~Stellato, G.~Banjac, P.~Goulart, A.~Bemporad, and S.~Boyd.
\newblock {OSQP}: an operator splitting solver for quadratic programs.
\newblock \emph{Mathematical Programming Computation}, 12\penalty0
  (4):\penalty0 637--672, 2020.
\newblock \doi{10.1007/s12532-020-00179-2}.
\newblock URL \url{https://doi.org/10.1007/s12532-020-00179-2}.

\bibitem[Tabuada(2009)]{tabuada2009verification}
Paulo Tabuada.
\newblock \emph{Verification and control of hybrid systems: a symbolic
  approach}.
\newblock Springer Science \& Business Media, 2009.

\end{thebibliography}

\end{document}